# ON THE UNIQUENESS OF THE INFINITE CLUSTER OF THE VACANT SET OF RANDOM INTERLACEMENTS

BY AUGUSTO TEIXEIRA

*ETH Zürich*

We consider the model of random interlacements on $\mathbb{Z}^d$ introduced in Sznitman [Vacant set of random interlacements and percolation (2007) preprint]. For this model, we prove the uniqueness of the infinite component of the vacant set. As a consequence, we derive the continuity in $u$ of the probability that the origin belongs to the infinite component of the vacant set at level $u$ in the supercritical phase $u < u_*$.

**1. Introduction.** In this paper we proceed with the study of percolation in the model of random interlacements introduced by Sznitman in [8]. This model is naturally linked to the analysis of random walks on the discrete torus $(\mathbb{Z}/N\mathbb{Z})^d$ ($d \geq 3$) and on the discrete cylinder $(\mathbb{Z}/N\mathbb{Z})^d \times \mathbb{Z}$ ($d \geq 2$), when the walk runs up to times proportional to $N^d$ and $N^{2d}$, respectively. Random interlacements heuristically describe the limiting microscopic "texture in the bulk" left by the random walk in these time scales, see [1, 3] and [9].

One can think of random interlacements as consisting of the trace of a cloud of trajectories in $\mathbb{Z}^d$ ($d \geq 3$) where a parameter $u \in (0, \infty)$ controls the density of the cloud. This pictorial description derives from the original construction given in [8] that will be detailed in Section 2.

The law $Q^u$ of the indicator function of the vacant set at level $u$, that is, the complement of the interlacement set, viewed as a random element of $\{0, 1\}^{\mathbb{Z}^d}$, is characterized as the unique law on $\{0, 1\}^{\mathbb{Z}^d}$, such that

(1.1) $\quad Q^u[Y_x = 1 \text{ for all } x \in K] = \exp\{-u \operatorname{cap}(K)\},$

for all finite sets $K \subset \mathbb{Z}^d$;

see [8], Sections 1 and 2. In the equation above, $(Y_x)_{x \in \mathbb{Z}^d}$ stands for the canonical coordinates on $\{0, 1\}^{\mathbb{Z}^d}$ (endowed with its canonical $\sigma$-algebra $\mathcal{Y}$) and $\mathrm{cap}(K)$ stands for the capacity of the set $K$; see (2.4). As shown in Theorem 2.1 of [8],

(1.2)     the translation operators are measure-preserving flows on $(\{0,1\}^{\mathbb{Z}^d}, \mathcal{Y}, Q^u)$ which are ergodic.

Our main interest here lies in the percolative properties of the random set $\{x \in \mathbb{Z}^d; Y_x = 1\}$. With this in mind, it is natural to introduce the function

$$\eta(u) = Q^u[0 \text{ belongs to an infinite component of } \{x \in \mathbb{Z}^d; Y_x = 1\}]$$

and the critical value

$$(1.3) \qquad u_* = \inf\{u \geq 0; \eta(u) = 0\}.$$

It was proved in [8] that for $d \geq 7$, $u_*$ is positive and that for $d \geq 3$, $u_*$ is finite; see Theorems 3.5 and 4.3. The positivity of $u_*$ in the case $d \geq 3$ was later established in [7], (0.5).

Let us mention here that the set $\{x \in \mathbb{Z}^d; Y_x = 0\}$ presents a different behavior. As shown in Corollary 2.3 of [8], which is proved using similar ideas to those of the present paper,

(1.4)     for every $u > 0$, $Q^u$-almost surely $\{x \in \mathbb{Z}^d; Y_x = 0\}$ is an infinite connected subset of $\mathbb{Z}^d$.

In this paper, we establish some properties of the critical-supercritical phases ($u \leq u_*$) of the vacant set. Our main result concerns the uniqueness of the infinite cluster.

THEOREM 1.1 ($d \geq 3$). *For arbitrary $u \geq 0$, the set $\{x \in \mathbb{Z}^d; Y_x = 1\}$ has $Q^u$-almost surely at most one infinite component.*

An important consequence of the theorem above is Corollary 1.2 regarding the continuity of the percolation probability in the supercritical regime.

COROLLARY 1.2 ($d \geq 3$). *The percolation probability $\eta(u)$ is a left-continuous function which is continuous on $[0, u_*)$.*

It is not known whether the vacant set percolates at the critical point $u_*$.

The proof of Theorem 1.1 roughly follows the strategy of Burton and Keane [2]. But their argument cannot be applied to the situation considered here, since it relies on the so-called finite energy property, that the measures $Q^u$ fail to satisfy

$$(1.5) \qquad 0 < Q^u(Y_y = 1 | Y_z, z \neq y) < 1, \qquad Q^u\text{-a.s., for all } y \in \mathbb{Z}^d;$$

for more details, see [4], Section 12, and [8], (1.60). Let us recall the two main steps in the argument of Burton and Keane. Using the fact that the number of infinite components is almost surely a constant (this is already contained in [6], Theorem 1), one first proves that this constant is not an integer bigger than 1 and finally that this constant is also finite. In both steps one makes use of the finite energy property to prescribe a specific configuration of occupied sites inside a box.

To give an intuition of the kind of restriction that we have in the present case, we recall from (1.4) that the set $\{x \in \mathbb{Z}^d; Y_x = 0\}$ has $Q^u$-a.s. no bounded component, so that we cannot arbitrarily assign values for the variables $Y_x$, not even for a finite set of them.

It is not clear how to overcome this difficulty with the description of the process that was given in terms of (1.1). Instead, we consider the characterization of the process by means of the interlacement set, that we pictorially described as a cloud of trajectories of paths in $\mathbb{Z}^d$. Loosely speaking, our scheme is to modify the behavior of these paths which eventually meet a previously chosen set $K \subset \mathbb{Z}^d$ in order to assemble the desired configuration inside $K$. This approach presents geometrical restrictions on the configurations we are able to design. Both to prove that the number of infinite components is finite and that it is not bigger than 1, we have to overcome these restrictions by carefully choosing the sets which we want to cover with the interlacement paths.

For the proof of the corollary, we adapt the argument of van den Berg and Keane [11].

This paper is organized as follows. In Section 2, we give the description of the model of random interlacement and quote some key results achieved in [8]. Section 3 contains the proofs of Theorem 1.1 and Corollary 1.2.

**2. A brief review of the random interlacements.** We start with some notation that will be used to describe the model of random interlacements. Let $\{e_k\}_{k=1,\ldots,d}$ denote the canonical basis vectors of $\mathbb{R}^d$ and $\pi_k : \mathbb{Z}^d \to \mathbb{Z}$ the projection on the $k$th direction. The euclidean norm is denoted by $|\cdot|$ and we say that two points $x$ and $y$ in $\mathbb{Z}^d$ are nearest neighbors $(x \leftrightarrow y)$ if $|x - y|$ equals 1. The closed ball $B(x,r) = \{y \in \mathbb{Z}; |y - x|_\infty \leq r\}$ is defined in terms of the $l_\infty$-norm $|x|_\infty = \sup\{|\pi_k(x)|; 1 \leq i \leq d\}$. For a set $K \subset \mathbb{Z}$, $|K|$ denotes its cardinality. We also distinguish the two boundaries of $K$:

(2.1)
$$\partial K = \{x \in K; x \text{ is a nearest neighbor of some point } y \in K^c\},$$
$$\partial_{\text{ext}} K = \{x \in K^c; x \text{ is a nearest neighbor of some point } y \in K\},$$

and define its "closure" as $\bar{K} = K \cup \partial_{\text{ext}} K$.

Throughout this paper the term path always refers to nearest-neighbor paths and by a piece of path we mean a finite sequence $(\tau_n)_{0 \leq n \leq l}$ of nearest-neighbor points. By entrance (resp. departure) points of a path $\tau : \mathbb{Z} \to \mathbb{Z}^d$ in a set $K \in \mathbb{Z}^d$, we mean

(2.2) $\qquad E = \{\tau(n); \tau(n-1) \in K^c \text{ and } \tau(n) \in K \text{ for } n \in \mathbb{Z}\}$

(2.3) $\qquad [\text{resp. } D = \{\tau(n); \tau(n+1) \in K^c \text{ and } \tau(n) \in K \text{ for } n \in \mathbb{Z}\}].$

The capacity of a finite set $K \subset \mathbb{Z}^d$ is given by the formula

$$\text{cap}(K) = \frac{1}{2d} \inf \left\{ \sum_{x \leftrightarrow y} |f(x) - f(y)|^2; f = 1 \text{ in } K, \right.$$

(2.4)

$$\left. f \text{ has finite support} \right\}.$$

Note above that, in the notation of weighted graphs, we assigned $1/2d$ for the conductance of the edges of $\mathbb{Z}^d$, so that the sum of the conductance of all edges incident to a fixed vertice is 1.

Consider the spaces $W_+$ and $W$ of infinite, respectively doubly infinite paths, tending to infinity at infinity in $\mathbb{Z}^d$ ($d \geq 3$):

$$W_+ = \left\{ \gamma : \mathbb{N} \to \mathbb{Z}^d; \gamma(i) \leftrightarrow \gamma(i+1) \text{ for each } i \geq 0, \right.$$

$$\left. \lim_{i \to \infty} |\gamma(i)| = \infty \right\},$$

(2.5)

$$W = \left\{ \gamma : \mathbb{Z} \to \mathbb{Z}^d; \gamma(i) \leftrightarrow \gamma(i+1) \text{ for each } i \in \mathbb{Z}, \right.$$

$$\left. \lim_{|i| \to \infty} |\gamma(i)| = \infty \right\},$$

and endow them with the $\sigma$-algebras $\mathcal{W}_+$ and $\mathcal{W}$ generated by the coordinate maps $X_n$. For a given path $w \in W_+$ (resp. $w' \in W$), we write $\text{Range}(w) = w(\mathbb{N})$ [resp. $\text{Range}(w') = w'(\mathbb{Z})$] and denote with $P_x$ the probability on $(W_+, \mathcal{W}_+)$ governing the simple random walk starting at $x$. Since $d \geq 3$, the simple random walk is transient, and the set $W_+$ supports the probability $P_x$. On $W_+$, the hitting time of a set $K \subset \mathbb{Z}^d$ is represented by

(2.6) $$\tilde{H}_K = \inf\{n \geq 1; X_n \in K\}.$$

Define the space $W^*$ of paths in $W$ modulo time-shift, that is,

(2.7) $W^* = W/\sim,$ where $\omega \sim \omega'$ if $\omega(\cdot) = \omega'(k + \cdot)$, for some $k \in \mathbb{Z}$,

and denote with $\pi^*$ the canonical projection from $W$ to $W^*$. The map $\pi^*$ induces a $\sigma$-algebra in $W^*$ given by $\mathcal{W}^* = \{A \subset W^*; (\pi^*)^{-1}(A) \in \mathcal{W}\}$.

The random interlacement is defined by means of a Poisson point process on $W^* \times \mathbb{R}_+$. The intensity of this process will be given by the product of a measure $\nu$ on $W^*$ and the Lebesgue measure; see, for instance, [10] for a reference about point processes. Here, the positive real line is introduced in order to couple different intensities of the point processes, allowing us to relate the probability of increasing events at different intensities $u, u' \in \mathbb{R}_+$ [an event $A \in \mathcal{Y}$ is said to be increasing if, whenever $\alpha \in A$ and $\alpha'(z) \geq \alpha(z)$ for all $z \in \mathbb{Z}^d$, we have $\alpha' \in A$]. The mentioned measure $\nu$ is constructed in [8] and in the next paragraph we give a characterization of it.

For any finite set $K$ of $\mathbb{Z}^d$, we write $W_K$ for the space of paths in $W$ that intersect $K$ and introduce $W_K^* = \pi^*(W_K)$. On $W_K$ we define the finite measure $Q_K$ such that, given $A$ and $B$ in the $\sigma$-algebra $\mathcal{W}_+$ and a piece of path $\tau : \{0, \ldots, l\} \to \mathbb{Z}^d$,

$$
\begin{aligned}
Q_K[(X_{-n})_{n \geq 0} &\in A, (X_n)_{0 \leq n \leq L_K} = \tau, (X_{n+L_K})_{n \geq 0} \in B] \\
&= P_{\tau(0)}[A, \tilde{H}_K = \infty] P_{\tau(0)}[(X_n)_{0 \leq n \leq l} = \tau] P_{\tau(l)}[B, \tilde{H}_K = \infty],
\end{aligned} \tag{2.8}
$$

where $L_K$ is the time of the last visit of a path to the set $K$.

A. S. Sznitman proved in [8] the consistency of the measures $\pi^* \circ Q_K$; that is, whenever $K \subset K'$, we have $1_{W_K^*} \cdot \pi^* \circ Q_{K'} = \pi^* \circ Q_K$. This readily implies the existence of a $\sigma$-finite measure $\nu$ in $W^*$ satisfying, for any finite set $K \subset \mathbb{Z}^d$,

$$1_{W_K^*} \cdot \nu = \pi^* \circ Q_K. \tag{2.9}$$

We are now in position to define the process we are interested in. Consider the space of point measures in $W^* \times \mathbb{R}_+$

$$
\Omega = \left\{ \omega = \sum_{i \geq 1} \delta_{(w_i^*, u_i)}; w_i^* \in W^*, u_i \in \mathbb{R}_+ \text{ and } \omega(W_K^* \times [0, u]) < \infty \right. \\
\left. \text{for every finite } K \subset \mathbb{Z}^d \text{ and } u \geq 0 \right\}, \tag{2.10}
$$

endowed with the $\sigma$-algebra generated by the evaluation maps $\omega \mapsto \omega(D)$ for $D \in \mathcal{W}^* \otimes \mathcal{B}(\mathbb{R}_+)$ [where $\mathcal{B}(\cdot)$ denotes the Borel $\sigma$-algebra]. Let $\mathbb{P}$ be the law of a Poisson point process in $W^* \times \mathbb{R}_+$ (see for this definition [10], Proposition 3.6) with intensity given by the product of $\nu$ with the Lebesgue measure on $\mathbb{R}_+$. For a given point measure $\omega$ in $\Omega$, written as $\omega = \sum_i \delta_{(w_i^*, u_i)}$, we define the *interlacement* and the *vacant set* at level $u$, respectively, as

$$\mathcal{I}^u(\omega) = \left\{ \bigcup_{i; u_i \leq u} \text{Range}(w_i^*) \right\}, \tag{2.11}$$

$$\mathcal{V}^u = \mathbb{Z}^d \setminus \mathcal{I}^u. \tag{2.12}$$

According to [8], Remark 2.2, in the notation of (1.1),

(2.13) $Q^u$ is the image measure of $(1_{[x \in \mathcal{V}^u]})_{x \in \mathbb{Z}^d}$ under $\mathbb{P}$.

So, the set $\{x \in \mathbb{Z}^d; Y_x = 1\}$ that was used to state Theorem 1.1 has the same law as the vacant set defined in (2.12).

We finish this section with a definition that will be crucial for the proofs in the next section.

DEFINITION 2.1. Given $u > 0$ and $\omega \in \Omega$, for a finite set $K \subset \mathbb{Z}^d$, we say that $x_1, \ldots, x_k \in K$ are an exiting $k$-tuple of $K$ at level $u$ if they belong to distinct infinite connected components of both $\mathcal{V}^u(\omega)$ and $\mathcal{V}^u(\omega) \setminus (K \setminus \{x_1, \ldots, x_k\})$.

**3. Uniqueness of the infinite component of the vacant set.** In this section we prove Theorem 1.1 and Corollary 1.2. In proving Theorem 1.1, just as for the case of independent Bernoulli percolation, we use the fact that the number of infinite components of $\mathcal{V}^u$ is almost surely a constant. We only need to discard the possibility that this constant is either an integer bigger than 1 or that it is infinite.

To handle the former possibility, the idea is to modify the behavior of the paths that visit a predefined box in order to join at least two of the infinite components touching it; see Lemma 3.2. This should not be done by simply removing all the paths that meet this box because this could change the configuration outside the box, and possibly create more infinite components of the vacant set. To handle the latter possibility, the rough idea is to build a trifurcation point (see Definition 3.3) using the interlacement paths; this is done in what follows Lemma 3.4

In both cases, the difficulty is due to the fact that the measures $Q^u$ fail to satisfy the so-called finite energy condition (1.5) used in [2]. For a more detailed discussion on this, we refer to (1.60) of [8].

We restate our main theorem in terms of $\mathcal{V}^u$.

THEOREM 3.1 ($d \geq 3$). *For any $u \geq 0$ with $\eta(u) > 0$, there exists $\mathbb{P}$-a.s. a unique infinite component in the vacant set $\mathcal{V}^u$.*

PROOF. Let $N^u$ be the number of infinite components of the set $\mathcal{V}^u$. Since this random variable is invariant under translations of the lattice, (2.13) and (1.2) imply that $N^u$ is $\mathbb{P}$-almost surely a constant $k(u)$. We start with:

LEMMA 3.2. *For any parameter $u$, $k(u) \in \{0, 1, \infty\}$.*

PROOF. We will write $k$ in place of $k(u)$ for simplicity. Suppose that for some $u$, we have $2 \leq k < \infty$. Then, by the continuity of $\mathbb{P}$, there exists an $L$ (that we suppose to be larger than 100) such that with positive probability the box $K = B(0, L)$ is connected to $k$ infinite components of $\mathcal{V}^u$.

To control the fact that the modifications we will do inside $\bar{K}$ [recall the definition below (2.1)] are not going to create more than $k$ infinite components in the vacant set, we first show that the event

(3.1) $A_1 = \{\bar{K} \cap \mathcal{I}^u \neq \varnothing$ and $\mathcal{V}^u \setminus \bar{K}$ has more than $k$ infinite components$\}$

has probability zero.

Define for this purpose a function $\phi_1 \colon W^*_{\bar K} \to W^*_{\bar K}$ in the following way. For a given path $w^* \in W^*_{\bar K}$, $\phi_1(w^*)$ is the path obtained by inserting in the moment of the first visit of $w^*$ in $\bar K$ a piece of path that covers the set $\bar K$.

With the help of (2.8) we see that $\phi_1 \circ (1_{W^*_{\bar K}} \cdot \nu)$ is absolutely continuous with respect to $1_{W^*_{\bar K}} \cdot \nu$. Extending $\phi_1$ as the identity on the complement of $W^*_{\bar K}$ and defining

$$(3.2) \quad \Phi_1(\omega) = \sum_{u_i < u} \delta_{(\phi_1(w^*_i), u_i)} + \sum_{u_i \geq u} \delta_{(w^*_i, u_i)}, \qquad \text{for } \omega = \sum_{i=1}^{\infty} \delta_{(w^*_i, u_i)},$$

we also have $\Phi_1 \circ \mathbb{P} \ll \mathbb{P}$.

For a given $\omega \in A_1$, the effect of applying $\Phi_1$ is to fill the set $K$, so that $A_1 \subset \Phi_1^{-1}[N^u > k]$. Then, since $N^u$ equals $k$ $\mathbb{P}$-almost surely, $\Phi_1 \circ \mathbb{P} \ll \mathbb{P}$ implies that $\mathbb{P}(A_1) = 0$.

Note that whenever $K$ is connected to $k$ infinite clusters at level $u$, it is possible to find in $\partial_{\text{ext}} K$ an exiting $k$-tuple of $\bar K$; recall the Definition 2.1. Since we have a finite number of choices for the points in this $k$-tuple, we can find nonrandom $z_1, \ldots, z_k \in \partial_{\text{ext}} K$ for which the event

$$(3.3) \qquad A_2 = \{z_1, \ldots, z_k \text{ is an exiting } k\text{-tuple of } \bar K \text{ at level } u\}$$

has positive probability.

One can choose a set $U \subset K \cup \{z_1, z_2\}$, containing a path joining $z_1$ and $z_2$ and such that $\mathcal{S} = \bar K \setminus U$ is connected (recall that we required $L$ to be larger than 100). For instance, take $U$ to be the union of: $K$ without its corners (i.e., $K \setminus \{x \in K; |x| = L\sqrt{d}\}$), the points $z_1$ and $z_2$ and the neighbors of $z_1$ and $z_2$ in $K$.

Given the choice of such $U$, we claim that it is possible to define a map $\phi_2 \colon W^*_{\bar K} \to W^*_{\bar K}$ such that

$$(3.4) \quad \begin{array}{l} \text{for any } w^* \text{ in } W^* \text{ for which all entrance and departure points} \\ \text{of } \bar K \text{ are in } \mathcal{S}, \phi_2(w^*) \text{ has the same range as } w^* \text{ out of } \bar K \text{ and} \\ \text{does not intersect } U. \end{array}$$

Indeed, if $w^*$ is as above, one can replace each of the excursions that $w^*$ performs in $\bar K$ by a piece of path contained in $\mathcal{S}$ with the same starting and ending points. For paths not satisfying the condition above, we make $\phi_2(w) = w$.

On the event $A_2$, both $z_1$ and $z_2$ are in the vacant set. Hence, given a point measure $\omega = \sum_{i \geq 1} \delta_{(w^*_i, u_i)} \in A_2$, the entrance and departure points of all paths $w^*_i$ in $\bar K$ are distinct from the $z_1$ and $z_2$ and therefore are in $\mathcal{S}$. As in the case of $\phi_1$, using (2.8) we can conclude that $\phi_2 \circ (1_{W^*_{\bar K}} \cdot \nu) \ll (1_{W^*_{\bar K}} \cdot \nu)$. If we extend $\phi_2$ as identity out of $W^*_{\bar K}$ and define $\Phi_2$ in analogy with (3.2) with $\phi_1$ replaced by $\phi_2$ we have also $\Phi_2 \circ \mathbb{P} \ll \mathbb{P}$.

Note that applying $\Phi_2$ to some $\omega \in A_2$, we join the vacant components associated with $z_1$ and $z_2$ using the set $U$ (that is now contained in $\mathcal{V}^u$) without changing the interlacements outside $\bar{K}$. Hence, in view of (3.1) and Definition 2.1, we conclude that $A_2 \setminus A_1 \subset \Phi_2^{-1}[N^u < k]$. So that

$$(3.5) \qquad 0 < \mathbb{P}[A_2 \setminus A_1] \leq \mathbb{P}[\Phi_2^{-1}[N^u < k]] = \Phi_2 \circ \mathbb{P}[N^u < k],$$

which, by the absolute continuity of $\Phi_2 \circ \mathbb{P}$ with respect to $\mathbb{P}$, implies that $\mathbb{P}[N^u < k] > 0$ and thus contradicts the assumption that $N^u$ is almost surely $k(u)$. Hence, $k \in \{0, 1, \infty\}$. □

Now we want to rule out the case $k = \infty$. We follow the broad strategy of Burton and Keane (see [2]), but in the absence of the finite energy property, a specific proof for the existence of trifurcation points needs to be developed.

DEFINITION 3.3. We call $y$ a trifurcation point at level $u$ if it belongs to an infinite component of $\mathcal{V}^u$ which is split into three distinct infinite components by the removal of $y$.

The rough strategy of the proof will be to modify the behavior of the interlacement paths to build corridors connecting three infinite clusters to a fixed point $y$.

For $L \geq 1$, denote the $(d-1)$-dimensional $l_\infty$-ball

$$(3.6) \qquad S_L = \{x \in \mathbb{Z}^d; \pi_1(x) = 0 \text{ and } |x|_\infty \leq L\},$$

see the beginning of Section 2 for the notation. In the next lemma we prove that for an appropriate choice of $L$, with positive probability we can find an exiting triple of $\bar{S}_L$ at level $u$ in which the tree points are mutually distant.

LEMMA 3.4. *Fix $u > 0$. Under the assumption that $k(u) = \infty$, there exists a positive integer $L_0$ and three points $z_1$, $z_2$ and $z_3$ in $\partial_{\text{ext}} S_{L_0}$, such that*

$$(3.7) \qquad |z_i - z_j| \geq 100d \qquad \text{for } i \neq j,$$

*and the event*

$$(3.8) \qquad A_3 = \{z_1, z_2 \text{ and } z_3 \text{ is an exiting triple of } \bar{S}_{L_0} \text{ at level } u\}$$

*has positive probability.*

PROOF. Given a connected component of $\mathbb{Z}^d$, we say that it is unbounded in the direction $+e_k$ (resp. $-e_k$) if its projection $\pi_k$ on the $k$th coordinate is unbounded from above (resp. from below). Clearly, every infinite component of $\mathcal{V}^u$ is unbounded in at least one of the $2d$ possible directions. Hence, when the number of unbounded components is equal to infinity, we can find one direction for which infinitely many components are unbounded, that is,

(3.9)
$$[N^u = \infty]$$
$$\subseteq \bigcup_{r \in \{+,-\}} \bigcup_{i=1}^{d} [\#\{\text{unbounded components in direction } re_i\} = \infty].$$

Since $[N^u = \infty]$ was supposed to be a $\mathbb{P}$-almost sure event, we can choose a fixed direction (without loss of generality we can assume it to be $+e_1$) for which the event in brackets in the right-hand side of the equation above has positive probability.

Denote with $F_l$ (for $l \in \mathbb{Z}$) the hyperplane perpendicular to $e_1$ containing the point $le_1$, $F_l = \{x \in \mathbb{Z}^d; \pi_1(x) = l\}$. It is clear that for every configuration in which there are infinitely many components unbounded in direction $+e_1$, we can take $l$ large enough such that $F_l$ intersects at least $M = 3|B(0, 100d)|$ of them, that is,

$$[\#\{\text{unbounded components in } +e_1\} = \infty]$$
$$\subseteq \bigcup_{l \geq 0} [\#\{\text{unbounded components in } +e_1 \text{ intersecting } F_l\} \geq M].$$

Hence, for a proper choice of $l_0$ (that we can assume to be 0 by the translation invariance of the process; see [8], Proposition 1.3), $F_{l_0}$ intersects at least $3|B(0, 100d)|$ infinite components of $\mathcal{V}^u$ with positive probability.

We can infer the same statement for some set $S_L \subset F_0$. Indeed, we note that as $L$ increases, the $(d-1)$-dimensional balls $S_L$ [recall the definition in (3.6)] eventually meet every component that touches $F_0$, so that

$$[\#\{\text{infinite components intersecting } F_0\} \geq M]$$
$$\subseteq \bigcup_{L \geq 0} [\#\{\text{inf. comp. intersecting } S_L\} \geq M]$$

and for a suitable choice of $L_0$,

(3.10)    $\mathbb{P}[\#\{\text{infinite components intersecting } S_{L_0}\} \geqslant M] > 0.$

On the event in brackets above, it is possible to find an exiting $M$-tuple of $\bar{S}_{L_0}$ at level $u$, and extract from it an exiting triple satisfying (3.7).

Since the number of possible choices for this triple is finite, we can find deterministic $z_1$, $z_2$ and $z_3$ in $\partial_{\text{ext}} S_{L_0}$ (with $|z_i - z_j| \geq 100d$ for $i \neq j$) that are an exiting triple of $\bar{S}_{L_0}$ at level $u$ with positive probability. $\square$

We now begin the construction of the trifurcation associated to this triple. Roughly speaking, we will use the interlacement paths to build tunnels connecting the $z_i$'s to some point $y$ in $S_{L_0}$ [recall from (3.6), that $S_{L_0}$ is a $(d-1)$-dimensional $l_\infty$-ball]. What we now describe are the paths $(\gamma_i)_{i=1,2,3}$ that are going to be the base for the construction of such tunnels.

Each $\gamma_i$ starts at the point $z_i$ ($\in \partial_{\text{ext}} S_{L_0}$). In the first step, the path $\gamma_i$ goes from $z_i$ to its unique neighbor in $S_{L_0}$. In the case this path did not reach the smaller ball $S_{L_0-1}$, it enters $S_{L_0-1}$ using as few steps as possible (at most $d-1$). Now, one can easily continue these paths from their entrance points in $S_{L_0-1}$ (that are still far from each other) to a common point $y$ without leaving $S_{L_0-1}$ and in such a way that $\mathcal{H} \setminus \{y\}$ has three connected components, where $\mathcal{H} = \bigcup_i \text{Range}(\gamma_i)$.

Our aim now is to fill the region $\mathcal{C} = \bar{S}_{L_0} \setminus \mathcal{H}$ that surrounds the $\gamma_i$'s using the paths of the interlacements. To do this, we first show that $\mathcal{C}$ is connected by partitioning it into the following four subsets:

$$\text{``top''} \ \mathcal{T} = \mathcal{C} \cap (S_{L_0} + e_1),$$
$$\text{``bottom''} \ \mathcal{B} = \mathcal{C} \cap (S_{L_0} - e_1),$$
$$\text{``kernel''} \ \mathcal{K} = \mathcal{C} \cap S_{L_0},$$
$$\text{``sides''} \ \mathcal{E} = \mathcal{C} \setminus (\mathcal{T} \cup \mathcal{B} \cup \mathcal{K}).$$

Since $\mathcal{T}$ and $\mathcal{B}$ are $(d-1)$-dimensional $l_\infty$-balls with at most three mutually distant points removed, they are connected. Every point of $\mathcal{K}$ is neighbor to both $\mathcal{T}$ and $\mathcal{B}$ and since $\mathcal{K}$ is nonempty, $\mathcal{T}$, $\mathcal{B}$ and $\mathcal{K}$ belong to the same component of $\mathcal{C}$. Finally, we show that all the components in the set $\mathcal{E}$ are connected in $\mathcal{C}$ to the kernel $\mathcal{K}$.

If $d \geq 4$, each of the $2d-2$ "faces" of $\mathcal{E}$ is a connected set just as $\mathcal{T}$ and $\mathcal{B}$ are. Indeed, these "faces" are $(d-2)$-dimensional $l_\infty$-balls with at most three far apart points removed. In this case, since in the construction of the paths $\gamma_i$ we allowed them to cover at most $3(d-1)$ points of $S_{L_0} \setminus S_{L_0-1}$, we have at least $(2L_0+1)^{(d-2)} - 3 - 3(d-1) > 0$ points in each of these faces that are neighbors to $\mathcal{K}$.

If $d = 3$, the set $\mathcal{E}$ is given by the union of four segments with at most three points removed. After this removal, each of them splits into at most four smaller segments. Consider one of these smaller segments which we denote by $P$. We show that $P$ is connected to $\mathcal{K}$ by considering two distinct cases. If the length of $P$ is greater than or equal to 100, we conclude as in the paragraph above that it has at least $100 - 3 \cdot (3-1) > 0$ points in the neighborhood of $\mathcal{K}$. If the length of $P$ is smaller than 100, with (3.7) we conclude that $P$ is a segment determined in one extreme by one of the points $z_i$'s and in the other extreme by the end of one of the four segments of length $2L_0 + 1$ that cover $\mathcal{E}$. By the particular construction we chose, the path associated with the point $z_i$ in the extreme of $P$ does not visit any of the neighbors of $P$ in $S_{L_0}$. Again by (3.7) we conclude that the other two points $z_i$ are far from $P$ and their respective paths $\gamma_i$'s spend at most their first $d-1$ steps out of $S_{L_0-1}$. Hence all the neighbors of $P$ in $S_{L_0}$ belong to $\mathcal{K}$.

Now we show how the connectedness of $\mathcal{C}$ proved above allows us to cover it with the interlacement paths. It implies, for instance, that for every pair of points $x$ and $y$ in $\partial_{\text{ext}} S_{L_0} \setminus \{z_1, z_2, z_3\}$, there exists a piece of path $\alpha_{x,y}$ connecting $x$ and $y$ and with $\text{Range}(\alpha_{x,y}) = \mathcal{C}$.

We claim that it is possible to define a map $\phi_3 : W^*_{\bar{S}_{L_0}} \to W^*_{\bar{S}_{L_0}}$ such that

(3.11) given a path $w^*$ for which all entrance and departure points in $\bar{S}_{L_0}$ belong to $\mathcal{C}$, $\phi_3(w^*)$ has the same range as $w^*$ out of $\bar{S}_{L_0}$ and $\text{Range}(\phi_3(w^*)) \cap \bar{S}_{L_0} = \mathcal{C}$.

Indeed, if $w^*$ satisfies the above property, one simply replaces all portions of the path performed inside $\bar{S}_{L_0}$ by the piece $\alpha$ (as above) corresponding to the entrance and departure points of this portion in $\bar{S}_{L_0}$. We set $\phi_3$ as the identity if $w^*$ does not satisfy the mentioned property.

We emphasize that the $w^*_i$ with $u_i \leq u$ appearing in a point measure $\omega \in A_3$ [see (2.10)] have the property that all its entrance and departure points of $\bar{S}_{L_0}$ are distinct from $z_1$, $z_2$ and $z_3$, and hence belong to $\mathcal{C}$.

Extend $\phi_3$ as the identity out of $W^*_{\bar{S}_{L_0}}$. Using (2.8) we conclude that $\phi_3 \circ \nu$ is absolutely continuous with respect to $\nu$. Defining $\Phi_3$ as in (3.2) with $\phi_1$ replaced by $\phi_3$, we note that also $\Phi_3 \circ \mathbb{P}$ is absolutely continuous with respect to $\mathbb{P}$ and on the event $\Phi_3(A_3)$, $y$ is a trifurcation point of $\mathcal{V}^u$.

Since $A_3$ has positive probability,

$$0 < \mathbb{P}(A_3) \leq \Phi_3 \circ \mathbb{P}[\Phi_3(A_3)],$$

implying that $\mathbb{P}(\Phi_3(A_3)) > 0$. The rest of the argument for uniqueness follows the proof in [2] with help of the ergodicity of translations mentioned in (1.2). This shows that $\mathbb{P}[N^u = \infty] = 0$, completing the proof of Theorem 1.1. □

An important consequence of the uniqueness of the infinite component is the continuity of the percolation probability in the supercritical phase, stated in Corollary 1.2. As mentioned in the Introduction, we follow the argument of van den Berg and Keane [11] with slight modifications.

PROOF OF COROLLARY 1.2. First we prove the left-continuity of $\eta$. The probability of the event $C^u_r = \{$the origin is connected by $\mathcal{V}^u$ to the boundary of the box $K = B(0, r)\}$ is a real analytic function of $u$ as follows, for instance, from the inclusion-exclusion formula (2.17) of [8].

Observe that in the notation of (1.3),

$$(3.12) \qquad \eta(u) = \mathbb{P}\left[\bigcap_{r \geq 1} C^u_r\right] = \lim_{r \to \infty} \mathbb{P}[C^u_r]$$

is a decreasing limit of continuous functions and hence is upper-semicontinuous on $\mathbb{R}_+$. Since it is monotone nonincreasing, $\eta$ is left-continuous on $[0, \infty)$.

To prove the right-continuity of $\eta$ in $[0, u_*)$, we need to understand the behavior of the event $C_\infty^u = \bigcap_{r \geq 1} C_r^u$ for $u$ in this interval. From now on, take a fixed $u < u_*$. By the continuity of $\mathbb{P}$ and the observation that $C_\infty^v$ is a nonincreasing family of events with respect to $v$, we have

$$(3.13) \qquad \lim_{v \downarrow u} \eta(v) = \lim_{v \downarrow u} \mathbb{P}[C_\infty^v] = \mathbb{P}\left[\bigcup_{v > u} C_\infty^v\right],$$

and all we need to prove is that this limit is in fact $\eta(u)$, or equivalently, that

$$(3.14) \qquad \mathbb{P}\left[C_\infty^u \setminus \left(\bigcup_{v > u} C_\infty^v\right)\right] = 0.$$

Fix a $v_o$ in $(u, u_*)$. Recalling the definition in (2.10), we consider the intersection of $C_\infty^u$ with

$$(3.15) \qquad B = \left\{\omega = \sum_i \delta_{(w_i^*, u_i)} \in \Omega; \text{ and there exist } J^u \text{ and } J^{v_o},\right.$$

$$\left. \text{unique infinite components of } \mathcal{V}^u \text{ and } \mathcal{V}^{v_o} \text{ respectively}\right\}.$$

In view of Theorem 1.1, $B$ is an almost sure event. As a step toward (3.14), we first show that

$$(3.16) \qquad B \cap C_\infty^u \subset \bigcup_{v > u} C_\infty^v.$$

By the uniqueness required in (3.15), for $\omega \in B$, $J^{v_o} \subset J^u$. Hence, whenever $\omega \in B \cap C_\infty^u$, we can find a finite path $\tau$ in $J^u$, connecting $J^{v_o}$ to the origin. It is clear by (2.10) that the number of pairs $(w_i^*, u_i)$ with $u_i < v_o$, such that $w_i^*$ touches $\tau$, is finite. Moreover, because $\text{Range}(\tau) \subset \mathcal{V}^u$, all the $u_i$'s in the pairs above are bigger than $u$, so that their infimum is bigger than $u$. Then, for parameters $v$ between $u$ and this infimum, the origin belongs to an infinite component. In other words, $\omega \in C_\infty^v$, so that $B \cap C_\infty^u \subset (\bigcup_{v > u} C_\infty^v)$ and (3.16) follows. The claim (3.14) is now a direct consequence of (3.16) and $\mathbb{P}[B] = 1$. □

REMARK 3.5. *(1) According to Remark 1.4 of [8], the construction of random interlacements can be straightforwardly generalized to other transient weighted graphs. It is a natural question whether the results of this paper can be achieved in some other situations. For instance, in what context*

UNIQUENESS IN RANDOM INTERLACEMENTS 13

*does Lemma 3.2 hold? Clearly the argument presented here is very specific to* $\mathbb{Z}^d$.

*(2) An important way to extend Theorem 1.1 is to give a quantitative analogous result for the supercritical phase. One could be interested, for example, in the decay with $n$, of the probability that inside $B(0,n)$ one can find two distinct components of diameter at least $\alpha \cdot n$, for $\alpha > 0$.*

*(3) One can also be interested in the question of simultaneous uniqueness. More precisely, is it true that $\mathbb{P}$-a.s. there exists at most one infinite component in the vacant set at level $u$ for every parameter $u > 0$? This is known to be true for Bernoulli independent percolation in several contexts; see, for instance, [5].*

**Acknowledgments.** We are grateful to Alain-Sol Sznitman, David Windisch and an anonymous referee for the valuable suggestions and careful reading.

Departement Mathematik
ETH Zürich HG G 10.1
Rämistrasse 101
8092 Zürich
Switzerland
E-mail: augusto.teixeira@math.ethz.ch